%% file: FreObs_article.tex
\documentclass[review,onefignum,onetabnum]{siamonline190516}

\input{FreObs_shared}

\begin{document}

\maketitle

\begin{abstract}
\input{sections/Abstract}
\end{abstract}



\section{Introduction} 
\label{sec:Introduction}
\input{sections/Introduction}

\section{Characterizing the Computational-Statistical Trade-Off} \label{sec:Characterization}
\input{sections/Characterization}

\section{Stochastic Approximation Methods for 4D-Var}
\label{sec:Stochastic}
\input{sections/Stochastic}

\section{Numerical Experiments}
\label{sec:Experiments}
\input{sections/Experiments}

\section{Conclusion}
\label{sec:Conclusion}
\input{sections/Conclusion}

\bibliographystyle{siamplain}
\bibliography{references}
\end{document}


\maketitle

%
%
%
%
%
%
%
%
%
%
%

\bibliographystyle{siamplain}
\bibliography{references}

%% file: FreObs_shared.tex
\usepackage{lipsum}
\usepackage{amsfonts}
\usepackage{graphicx}
\usepackage{subcaption}
\usepackage{epstopdf}
\usepackage{algorithmic}
\usepackage{cite}
\usepackage{tikz}
\usepackage{pgfplots}
\usepackage{pgfplotstable}
\usepackage{cleveref}
\usetikzlibrary{decorations.pathreplacing}
\usepackage{csvsimple}
\usepackage{siunitx}
\usepackage{array, booktabs, makecell, caption}
\usepackage[figuresright]{rotating}
\ifpdf
\DeclareGraphicsExtensions{.eps,.pdf,.png,.jpg}
\else
\DeclareGraphicsExtensions{.eps}
\fi

\usepackage{enumitem}
\setlist[enumerate]{leftmargin=.5in}
\setlist[itemize]{leftmargin=.5in}

\newcommand\myhead[1]{\parbox[t]{5.8em}{\centering\bfseries#1\par\kern3mm}}

\newsiamremark{remark}{Remark}
\newsiamremark{example}{Example}

\newsiamthm{assumption}{Assumption}
\crefname{assumption}{Assumption}{Assumptions}

\headers{High-Frequency Observations in Data Assimilation}{Zhang \& Patel}

\title{Stochastic Approximation for High-frequency Observations in Data Assimilation\thanks{Submitted to the editors on November 4, 2020.\funding{Vivak Patel is supported by the Wisconsin Alumni Research Foundation.}}}

\author{Shushu Zhang\thanks{Department of Statistics, University of Wisconsin, Madison, WI
  (\email{szhang695@wisc.edu}).}
\and Vivak Patel\thanks{Department of Statistics, University of Wisconsin, Madison, WI 
  (\email{vivak.patel@wisc.edu}, \url{vivakpatel.org}).}}

\usepackage{amsopn}

\DeclareMathOperator{\st}{subject~to}

\newcommand{\norm}[1]{\left\Vert #1 \right\Vert}
\newcommand{\E}[1]{\mathbb{E}\left[ #1 \right]}
\newcommand{\Prb}[1]{\mathbb{P}\left[ #1 \right]}
\newcommand{\1}[1]{\mathbf{1}\left[ #1 \right]}

%% file: sections/Abstract.tex
With the increasing penetration of high-frequency sensors across a number of biological and physical systems, the abundance of the resulting observations offers opportunities for higher statistical accuracy of down-stream estimates, but their frequency results in a plethora of computational problems in data assimilation tasks. The high-frequency of these observations has been traditionally dealt with by using data modification strategies such as accumulation, averaging, and sampling. However, these data modification strategies will reduce the quality of the estimates, which may be untenable for many systems. Therefore, to ensure high-quality estimates, we adapt stochastic approximation methods to address the unique challenges of high-frequency observations in data assimilation. As a result, we are able to produce estimates that leverage all of the observations in a manner that avoids the aforementioned computational problems and preserves the statistical accuracy of the estimates.


%% file: sections/Introduction.tex
Data assimilation (DA) has become a corner stone for estimating time-varying phenomena owing to advancements in programming paradigms, numerical methods, differential equation modeling, and sensor technology \cite{simon2006optimal}. 
Owing to the high expense of numerical integration, current, state-of-the-art DA developments focus primarily on advancing both code bases and techniques for the high-fidelity integration of the differential equation models (e.g., \cite{ecmwf2017four,noaa2016four}). 
While such high-fidelity integration allows for better accuracy to the underlying model, the overall DA solution accuracy is also fundamentally limited by the DA problem's implicit statistical accuracy; 
that is, the overall DA problem's solution accuracy can also be limited by the number of observations available or \textit{by how the observations are incorporated into the DA problem formulation}. 

With increased penetration of high-frequency sensor technology across many application areas (e.g., \cite{langland2004estimation,sutton2014high,monti2016international}), the DA problem's statistical accuracy is unencumbered by the number of observations available. 
On the other hand, owing to the frequency of such sensors, the DA problem is subject to \textit{high-frequency observations}: those observations that occur at time points shorter than a desirable numerical integration step size. As a result, solving the resulting DA problem would often exceed a desirable computational budget even with such techniques as checkpointing and windowing \cite{griewank2012invented}.
Therefore, to control computational costs, the DA problem is often altered by applying such data-modification techniques as sampling\footnote{``Sampling" is also referred to ``thinning" and ``subsetting" in different literature.} the data (e.g., \cite{sutherland2014weather}), or displacing the data to more convenient spatio-temporal locations (e.g., \cite{courtier1998ecmwf,rabier1998ecmwf}).
Consequently, the DA solution's statistical accuracy is compromised as these data modifications either introduce bias to, or increase the variance of, the solution (see Ch. 8 \& Ch. 9 of \cite{ljung1999system}).

In summary, DA solutions may have greater accuracy to the underlying model owing to improvements in numerical integration, but the solution's overall statistical accuracy is compromised owing to modifications of the data so as to control computational costs. 
Especially in the high-frequency observation regime where there is a need to leverage modern high-frequency sensor technologies to develop more granular insights (e.g., \cite{langland2004estimation,sutton2014high,monti2016international}), compromising the DA solution's statistical accuracy to control the computational burden would hinder such insights. Therefore, in the high-frequency observation regime, this DA problem's statistical-computational trade-off raises several practical questions: 
\begin{remunerate}
\item How severe is the trade-off between statistical accuracy and computational expediency under different, common data modification (e.g., sampling, discarding, displacing) strategies?
\item Is this statistical-computational trade-off a fundamental feature of the computational nature of the DA problem or is it possible to circumvent this trade-off?
\end{remunerate}

The first question has been partially investigated numerically within the numerical weather prediction community with a focus on observation errors, correlation treatment, and specific data modification strategies \cite{bergman1976analysis,liu2002interaction,liu2003potential,hoffman2018effect}, but has not be analyzed under all of the common data modification strategies. Therefore, here, we add to this literature by supplying additional characterizations of the trade-off between statistical accuracy and computational burden for the DA problem under the most common data modification strategies. As a result, we will be able to provide some additional guidance to practitioners who need to understand the consequences of their actions when they navigate the aforementioned trade-off. Additionally and more importantly, to make progress towards addressing the second question, we develop methods that can maintain a competitive computational budget, while also ensuring that statistical accuracy is not compromised. Specifically, we extend recent developments in stochastic optimization for large data sets (e.g., \cite{bottou2018optimization}) to the more complex setting of the DA problem with high-frequency observations. As a result, we will supply practitioners a new set of tools to navigate the statistical-computational trade-off. To summarize,
\begin{remunerate}
\item We supply additional numerical evidence for the computational-statistical trade-off in DA with a more comprehensive investigation of common data modification schemes; and
\item We extend stochastic approximation methods to the more complex setting of Data Assimilation with high-frequency observations.
\end{remunerate}

The remainder of this paper is organized as follows. In \cref{sec:Characterization}, we characterize the statistical-computational trade-off for common data modification methods. In \cref{sec:Stochastic}, we introduce stochastic approximation methods to circumvent the statistical-computational trade-off. In \cref{sec:Experiments}, we compare traditional solvers with and without data modification against the methods that we introduce in \cref{sec:Stochastic}. In \cref{sec:Conclusion}, we conclude this work.

%% file: sections/Characterization.tex
To control the computational costs of solving a data assimilation (DA) problem with high-frequency observations---that is, those observations that occur at time scales shorter than the desired numerical integration time step---, data modification strategies such as sampling or displacing observations are employed, which then compromise the statistical accuracy of the solution to the DA problem. While such a trade-off between the computational costs of solving the problem and the resulting statistical accuracy of the problem are not ideal, data modification strategies may or may not result in a tolerable loss of statistical accuracy (at scales at which the desired quantity is being estimated) \cite{bergman1976analysis}. Therefore, here, we will supply additional numerical-evidence-based guidance to understand this computational-statistical trade-off for common data modification strategies, which builds on the more focused work from the numerical weather prediction literature \cite{bergman1976analysis,liu2002interaction,liu2003potential,hoffman2018effect}.

We begin by describing the DA problem in the current context and these common data modification strategies. We then study a simple system that provides basic guidance on how these data modifications can impact the accuracy of the solution.

\subsection{DA Problem with High Frequency Observations}
To understand the consequences of the data modification schemes on the data assimilation problem with high-frequency observations, we must first be clear about what we mean by the data assimilation problem under high-frequency observations. Broadly, the data assimilation problem is a maximum likelihood estimation task or a maximum posterior distribution estimation task (depending on whether one is a Frequentist or Bayesian) with a dynamical system constraint, where the quantity being estimated is determined by the end goal of the data assimilation problem (i.e., filtering, smoothing, predicting). 

Note, the statistically rigorous data assimilation problem as described below (see \cref{eqn-optimization-problem}) is often referred to as the variational approach in the state estimation literature in order to distinguish it from approximate filtering methods (e.g., Extended Kalman filter, Gaussian Ensemble Filters, Particle Filters). However, in the high-frequency regime, these filtering methods ought to be viewed as solutions to the variational data assimilation problem subject to a data modification strategy under a circumscribed assimilation window. For example, the Kalman filter applied to an assimilation problem with high-frequency observations is often modified by shifting the observations to the next integration time point, which is equivalent to stating the full variational problem where the integration step size is the assimilation window and the accumulation data modification strategy (described presently) is employed.

\subsection{Common Data Modification Strategies} \label{subsection-data-modification}
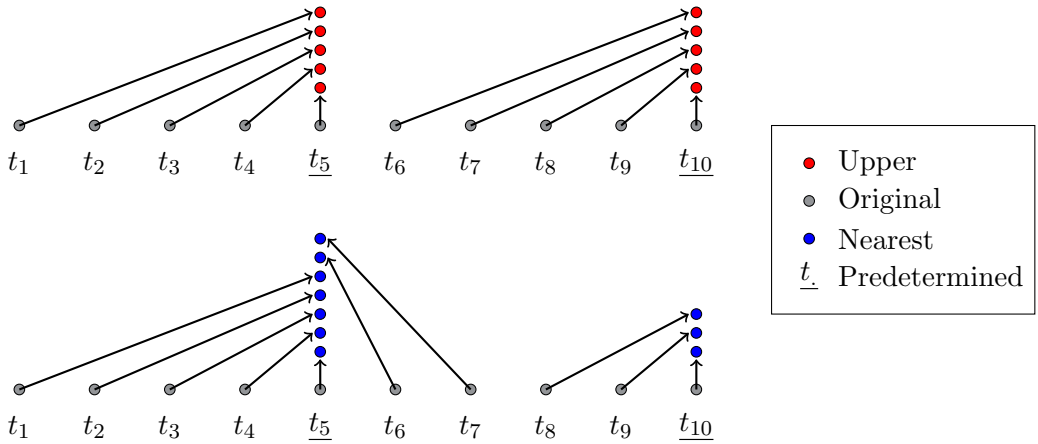
\begin{figure}[ht]
\input{tikz/accumulation}
\caption{In the upper graphic, the ``upper'' accumulation strategy is shown: in this approach, the original observation times (gray) between two adjacent predetermined times (with underline) are displaced to the ``upper'' time points (red) of the interval.
In the lower graphic, the ``nearest'' accumulation strategy is shown: in this approach, the original observation times (gray) are displaced to the ``nearest'' time point (blue) from the predetermined time sets (with underline).}
\label{figure:accumulation}
\end{figure}

Data modification strategies fall under three general categories: accumulation, averaging, and sampling.
The accumulation data modification approach simply displacing the observation time of a given data point to a more convenient location as shown in \cref{figure:accumulation}.\footnote{We can also talk about shifting the geospatial location of the observations \cite{salonen2009doppler,xu2011measuring}, but this is not as relevant in the context of high frequency observations.  In general, geospatial modifications are far more common than chronological modifications.} Note, the accumulation data modification approach does not modify the observed values themselves, it only shifts the observation times.

Unlike the accumulation strategy, the averaging strategy selects a set of observations within some time interval, computes their average (if possible), and assigns the value to a time point in the specified interval. The averaging strategy is more commonly observed in practice \cite{salonen2009doppler,xu2011measuring}, and often appears under the name of superobservations \cite{hoffman2018effect}. The averaging strategy is pictorially represented in \cref{figure:averaging}.

\begin{figure}[hp]
\begin{minipage}{1.0\textwidth}
\input{tikz/averaging}
\caption{In the upper graphic, the ``upper'' averaging strategy is shown: the original values (gray) between two adjacent predetermined thinned time points (with underline) are averaged (red) and assigned to the upper time point of the interval. In the lower graphic, the ``nearest'' averaging strategy is shown: the original observation values (gray) nearest some predetermined time point (with underline) are averaged (blue) and assigned to this time point. Note, the radius of the circles is representative of the number of data points used in the averaging scheme.}
\label{figure:averaging}
\end{minipage}

\begin{minipage}{1.0\textwidth}
\input{tikz/subsetting}
\caption{In the upper graphic, a systematic strategy is shown: original observation (gray) is selected by a predetermined constant distance regardless of the quality of the observation---the selected subset of observations is shown in red. In the lower graphic, an adaptive strategy is shown: the original observations (gray) with greater information content (e.g., higher precision), as shown by the larger radius of the circle, are given preference in the selection process, which results in the selected subset shown in blue.}
\label{figure:sampling}
\end{minipage}
\end{figure}

The sampling strategy is the most popular data modification approach, and comes in several varieties in the literature. The systematic sampling approach only keeps regularly spaced observations and discards the remaining observations \cite{li2010model,liu2002interaction,liu2003potential,petersen1963representative,singh2013practical,sutherland2014weather}. As a more recent development, adaptive sampling attempts to leverage some measure of information content to choose which observations to keep and which ones to discard \cite{gratton2015observation,ochotta2005adaptive,reale2018impact}. These adaptive sampling methods appear to be closely linked to matrix sketching that has become popular in linear algebra \cite{woodruff2014sketching}, but we will not comment on this further here. The sampling strategy is pictorially represented in \cref{figure:sampling}.

\subsection{Consequences of Data Modification} \label{subsection-data-modification-consequences}
In order to understand the performance of the aforementioned data modification schemes in terms of loss in accuracy and gain in computation budget, we supply numerical evidence on a simple example of a relaxation oscillator: the FitzHugh-Nagumo model for neuronal excitability as characterized in \cref{eq:FN}. 
\begin{equation}\label{eq:FN}
\begin{aligned}
\dot{v} &= v - v^3/3 - w + ii\\
\dot{w} &= (v - a - bw)/\tau,
\end{aligned}
\end{equation}
where $v$ is the membrane potential; $w$ is a recovery variable; and a four-dimensional unknown parameter variable $ii, a, b, \tau$ \cite{rocsoreanu2012fitzhugh}, which we integrate over a time interval of fifty seconds. 

To supply this numerical evidence, we can choose two routes. In the first route, we can \textit{ignore} the computational costs of determining the maximum likelihood estimator for different data modification schemes, and state what the consequences are of the data modification schemes. In the second route, we can fix a computational budget (using time as a proxy), allow the procedures to run, and then state the consequences of the different procedures. In the latter route, we are subject to the procedures that we choose and the implementation, whereas, in the former route, we do not suffer from having to make this choice. Therefore, for now, we take the former route which is method independent, and we save the latter approach to the experiments below.

To perform the experiments, we choose a ``true'' parameter value; we simulate the dynamics of the system; we generate high-frequency observations (one hundred per second) from linear transformations of the resulting states of the system plus independent Gaussian noise; and, finally, we solve the resulting maximum likelihood problem (initialized at the ``true'' parameter) with and without the aforementioned data modification strategies. In order to compare the resulting parameter estimates, we report the values of
\begin{equation} \label{eqn-relative-accuracy}
\frac{G_{\mathrm{no.mod}}(\hat \theta_{\mathrm{mod}}) - G_{\mathrm{no.mod}}(\hat \theta_{\mathrm{no.mod}}) }{G_{\mathrm{no.mod}}(\hat \theta_{\mathrm{no.mod}})},
\end{equation}
where $G_{\mathrm{no.mod}}(\theta)$ is the negative loglikelihood of the stated problem without any data modification; $\hat \theta_{\mathrm{no.mod}}$ is the maximum likelihood estimator from optimizing the likelihood without data modification (i.e., minimizing $G_{\mathrm{no.mod}}(\theta)$); and $\hat \theta_{\mathrm{mod}}$ is a placeholder for a parameter estimate resulting from minimizing the negative loglikelihood with data modification. Moreover, we also report the percentage of observation time points (POTP)---either 1.0\% or 10\%---used for the data modification schemes relative to the number of observation time points in original set of observations. These values are reported in \Cref{tab:observation error}. 

\begin{table}[hbt]
	\caption{Relative Error by Data Modification Scheme}
	\label{tab:observation error}
	\centering 
	\sisetup{scientific-notation = true, round-mode=places,round-precision=2}
	\begin{filecontents*}{"data/FitzHugh_Nagumo_Relative_Error.csv"}
x1,x2,x3,x4,x5,x6,x7
1020.3493075740903,1.1603505080168057,0.8949674523997444,1.5560745195811432,0.636211617725459,0.0004629782444496368,1.0
0.021386800155360035,4.1497682002667515e-5,0.009096064199293919,5.4026358853253205e-5,0.0004629782444035749,0.0003621088369267077,10.0
\end{filecontents*}
\csvreader[
	column count=4,
	tabular={p{0.5in}p{1in}p{1in}p{1in}}, 
	table head=\toprule  \textbf{POTP (\%)} & \textbf{Accumulate Upper} & \textbf{Accumulate Nearest} & \textbf{Average \mbox{Upper}} \\\midrule,
	table foot=\bottomrule	
	]{"data/FitzHugh_Nagumo_Relative_Error.csv"}
	{}
	{\csvcolvii & \num{\csvcoli} & \num{\csvcolii} & \num{\csvcoliii}} \\\bigskip 

	\csvreader[
	column count=4,
	tabular={p{0.5in}p{1in}p{1in}p{1in}}, 
	table head=\toprule  \textbf{POTP (\%)} & \textbf{Average Nearest} & \textbf{Simple ~ Random Sampling} & \textbf{Systematic Random Sampling} \\\midrule,
	table foot=\bottomrule	
	]{"data/FitzHugh_Nagumo_Relative_Error.csv"}
	{}
	{\csvcolvii & \num{\csvcoliv} & \num{\csvcolv} & \num{\csvcolvi}}
\end{table}

As summarized in \cref{tab:observation error}, each of these data modification schemes comes with a reduction in the quality of the approximation, though to varying degrees. When we are sampling the data at one percent, which roughly corresponds to the ideal integration time step of one second, systematic random sampling produces a negligible relative error, while the next best relative error is given by the simple random sampling scheme at 64\%. Of course, both of these random sampling schemes are highly subject to variability, and the randomness of both schemes can also produce other side effects: as we describe below, simple random sampling is a rather expensive scheme for data assimilation tasks; and \textit{a single} systematic random sample can miss periodic behaviors occurring at smaller time scales. Indeed, we need to increase the computational expense by a (generally untenable) factor of ten to see more reasonable relative errors for the data modification schemes and to stave off the effects of the randomness in the random sampling schemes. In summary, we see that most of the data modification schemes perform rather poorly or variably at the preferred numerical integration step size, and require a substantial and perhaps undermining amount of additional computational resources for reasonable and reliable performance.

%% file: tikz/accumulation.tex
\begin{tikzpicture}[level/.style={},decoration={brace,mirror,amplitude=7}]	
	\node[right] at (-0.5,2.5) {\bf Data Modification by Accumulation};	
	
	\draw (11,-2.5) rectangle (14.5,0);
	
		\draw[fill=red] (11.5,-0.5) circle (2pt);
		\node[right] at (11.75,-0.5) {Upper};
		
		\draw[fill=gray] (11.5,-1) circle (2pt);
		\node[right] at (11.75,-1) {Original};

		\draw[fill=blue] (11.5,-1.5) circle (2pt);
		\node[right] at (11.75,-1.5) {Nearest};
		
		\node at (11.5,-2) {\underline{$t_{.}$}};
		\node[right] at (11.75,-2) {Predetermined};
	
		
	\foreach \x in {1,2,3,4,5,6,7,8,9,10} 
		\draw[fill=gray] (\x,0) circle (2pt);
	\foreach \x in {1,2,3,4,6,7,8,9} 
		\node at (\x,-0.5) {$t_{\x}$};	
	\foreach \x in {5,10}
		\node at (\x,-0.5) {\underline{$t_{\x}$}};	

	\foreach \y in {0.5,0.75,1.0,1.25,1.5}
		\draw[fill=red] (5,\y) circle (2pt);
	\foreach \y in {0.5,0.75,1.0,1.25,1.5}
		\draw[fill=red] (10,\y) circle (2pt);
	
	\foreach \x/\y in {1/1.5,2/1.25,3/1.0,4/0.75}
		\draw [thick,->] (\x,0) -- (4.9,\y);
	\draw [thick,->] (5,0) -- (5,0.4);
	
	\foreach \x/\y in {6/1.5,7/1.25,8/1.0,9/0.75}
		\draw [thick,->] (\x,0) -- (9.9,\y);
	\draw [thick,->] (10,0) -- (10,0.4);
	
	
	\foreach \x in {1,2,3,4,5,6,7,8,9,10} 
		\draw[fill=gray] (\x,-3.5) circle (2pt);
	\foreach \x in {1,2,3,4,6,7,8,9} 
		\node at (\x,-4) {$t_{\x}$};
	\foreach \x in {5,10}
		\node at (\x,-4) {\underline{$t_{\x}$}};	

	\foreach \y in {0.5,0.75,1.0,1.25,1.5,1.75,2.0}
		\draw[fill=blue] (5,-3.5 + \y) circle (2pt);
	\foreach \y in {0.5,0.75,1.0}
		\draw[fill=blue] (10,-3.5 + \y) circle (2pt);
		
	\foreach \x/\y in {1/1.5,2/1.25,3/1.0,4/0.75}
		\draw [thick,->] (\x,-3.5 + 0) -- (4.9,-3.5 + \y);
	\draw [thick,->] (5,-3.5 + 0) -- (5, -3.5 + 0.4);
	\foreach \x/\y in {6/1.75,7/2.0}
		\draw [thick,->] (\x,-3.5 + 0) -- (5.1,-3.5 + \y);
	
	\foreach \x/\y in {8/1.0,9/0.75}
		\draw [thick,->] (\x,-3.5 + 0) -- (9.9,-3.5 + \y);
	\draw [thick,->] (10,-3.5 + 0) -- (10,-3.5 + 0.4);
	
\end{tikzpicture}

%% file: tikz/averaging.tex
\begin{tikzpicture}[level/.style={},decoration={brace,mirror,amplitude=7}]	
	\node[right] at (-0.5,2.5) {\bf Data Modification by Averaging};	
	
	\draw (11,-2.5) rectangle (14.5,0);
	
		\draw[fill=red] (11.5,-0.5) circle (2pt);
		\node[right] at (11.75,-0.5) {Upper};
		
		\draw[fill=gray] (11.5,-1) circle (2pt);
		\node[right] at (11.75,-1) {Original};
		
        \draw[fill=blue] (11.5,-1.5) circle (2pt);
        \node[right] at (11.75,-1.5) {Nearest};

        \node at (11.5,-2) {\underline{$t_{.}$}};
        \node[right] at (11.75,-2) {Predetermined};
	
		
	\foreach \x in {1,2,3,4,5,6,7,8,9,10} 
		\draw[fill=gray] (\x,0) circle (2pt);
	\foreach \x in {1,2,3,4,6,7,8,9} 
		\node at (\x,-0.5) {$t_{\x}$};
	\foreach \x in {5,10}
	    \node at (\x,-0.5) {\underline{$t_{\x}$}};	
	
	\draw[fill=red] (5,1.0) circle (4pt);
	\draw[fill=red] (10,1.0) circle (4pt);
	
	\foreach \x in {1,2,3,4}
		\draw [thick,->] (\x,0) .. controls (\x + 0.5,1) .. (4.8,1.0);
	\draw [thick,->] (5,0) -- (5,0.8);
	
	\foreach \x in {6,7,8,9}
		\draw [thick,->] (\x,0) .. controls (\x + 0.5,1) .. (9.8,1.0);
	\draw [thick,->] (10,0) -- (10,0.8);
	
	
	\foreach \x in {1,2,3,4,5,6,7,8,9,10} 
		\draw[fill=gray] (\x,-2.5) circle (2pt);
	\foreach \x in {1,2,3,4,6,7,8,9} 
		\node at (\x,-3) {$t_{\x}$};
	\foreach \x in {5,10}
	\node at (\x,-3) {\underline{$t_{\x}$}};	
	
	\draw[fill=blue] (5,-1.5) circle (7pt);
	\draw[fill=blue] (10,-1.5) circle (3pt);
		
	\foreach \x in {1,2,3,4}
		\draw [thick,->] (\x,-2.5) ..controls (\x + 0.5,-1.5) .. (4.7,-1.5);
		
	\draw [thick, ->] (5, -2.5) -- (5,-1.8);

	\foreach \x in {6,7}
		\draw [thick,->] (\x,-2.5) ..controls (\x - 0.5,-1.5) .. (5.3,-1.5);
	
	\foreach \x in {8,9}
		\draw [thick,->] (\x,-2.5) ..controls (\x + 0.5, -1.5) .. (9.8,-1.5);

	\draw [thick,->] (10,-2.5) -- (10,-1.7);
	
\end{tikzpicture}

%% file: tikz/subsetting.tex
\begin{tikzpicture}[level/.style={},decoration={brace,mirror,amplitude=7}]	
	\node[right] at (-0.5,2) {\bf Data Modification by Sampling};	
	
	\draw (11,-2) rectangle (14,0);
	
		\draw[fill=red] (11.5,-0.5) circle (2pt);
		\node[right] at (11.75,-0.5) {Systematic};
		
		\draw[fill=gray] (11.5,-1) circle (2pt);
		\node[right] at (11.75,-1) {Original};
		
		\draw[fill=blue] (11.5,-1.5) circle (2pt);
		\node[right] at (11.75,-1.5) {Adaptive};
	
		
	\foreach \x in {1,2,3,4,5,6,7,8,9,10} 
		\draw[fill=gray] (\x,0) circle (2pt);
	\foreach \x in {1,2,3,4,5,6,7,8,9,10} 
		\node at (\x,-0.5) {$t_{\x}$};
	
	\foreach \x in {1,3,5,7,9}
		\draw[fill=red] (\x,1.0) circle (2pt);
	
	\foreach \x in {1,3,5,7,9}
		\draw [thick,->] (\x,0) -- (\x,0.9);
	
	
	\foreach \x in {1,4,5,6,10} 
		\draw[fill=gray] (\x,-2.5) circle (2pt);
	\foreach \x in {2,3,7,8}
		\draw[fill=gray] (\x,-2.5) circle (4pt);
	\foreach \x in {9}
		\draw[fill=gray] (\x,-2.5) circle (6pt);
	\foreach \x in {1,2,3,4,5,6,7,8,9,10} 
		\node at (\x,-3) {$t_{\x}$};

	\foreach \x in {1}
		\draw[fill=blue] (\x,-1.5) circle (2pt);
	\foreach \x in {2,3,7}
		\draw[fill=blue] (\x,-1.5) circle (4pt);
	\foreach \x in {9}
		\draw[fill=blue] (\x,-1.5) circle (6pt);	
			
	\foreach \x in {1,2,3,7,9}
		\draw [thick,->] (\x,-2.5) -- (\x,-1.8);
	
\end{tikzpicture}

%% file: sections/Stochastic.tex
Given that we have underscored the challenges with data modification methods, we now propose a methodology, based on stochastic gradients, to address the shortfalls from data modification methods. In particular, we adapt a first-order stochastic gradient methodology and a second-order stochastic gradient methodology. For the first-order method, we adapt Stochastic Gradient Descent (SGD) \cite{robbins1951stochastic,blum1954multidimensional,chung1954stochastic}; for the second-order method, we adapt a stochastic method that correspond to Gauss-Newton methods \cite{bertsekas1996incremental,amari1998natural,patel2016kalman}, which we refer to as Kalman-based Stochastic Gradient Descent (kSGD). Note, we consider this particular second-order method owing to its ability to exploit common structures that appear in data assimilation problems.

It is important to recognize that we are not simply applying these stochastic methods directly to the data assimilation problem: a direct application would, as we explain below, be a rather inefficient use of computational resources. To this end, we begin by formulating the precise set of problems that we are interested in analyzing; we then discuss the methodology for computing the stochastic gradients for SGD; and, finally, we discuss the methodology for computing the stochastic gradients for kSGD.

\subsection{Problem Formulation}

We are interested in determining a parameter $\theta \in \mathbb{R}^p$ that governs a collection of variables $\lbrace x(t) : t \in [0,T] \rbrace \subset \mathbb{R}^d$ that satisfy
\begin{equation}
\dot x(t) = f(t, x(t), \theta),
\end{equation}
and some initial condition that may depend on $\theta$,\footnote{We can allow the dependence to be any smooth function of $\theta$. We will not add this additional complexity here, but it is readily addressable.} where $f$ is a known function that we will further specify below. Without loss of generality, we can augment the state vector $x$ by the parameter $\theta$, augment the dynamics $f$ by the $\mathbb{R}^p$-zero vector, 
and recast the problem as determining the initial value of the augmented state. To be specific, with an abuse of notation, we can rewrite the dynamics as
\begin{equation} \label{eqn-dynamics-ivp}
\left\lbrace
\begin{aligned}
\dot x(t) &= f(t, x(t)) \\
x(0) &= \theta,
\end{aligned} \right.
\end{equation}
where $\theta \in \mathbb{R}^{q}$ and $q = p + d$.

Now, to determine $\theta$ in \cref{eqn-dynamics-ivp}, we will have access to a set of possibly noisy observations $\lbrace y_i \rbrace \subset \mathbb{R}^n$ that are measured at an increasing sequence of time points $t_i \in [0,T]$, respectively. Moreover, we assume that the observations are each dependent only on the state at the corresponding observation time. In order to incorporate these observations, we assume that there is some well-selected loss function that relates the observations to the corresponding states, which we denote by $\ell: \mathbb{R}^{n} \times \mathbb{R}^q \to \mathbb{R}$. 
With these components, we can state the estimation problem as
\begin{equation} \label{eqn-optimization-problem}
\begin{aligned}
\min_\theta & ~ G(\theta) := \sum_{i=1}^N \ell( y_i, x(t_i) ) \\
\st & ~ \dot x (t) = f(t, x(t)), ~ t \in [0,T] \\
			 & ~ x(0) = \theta.
\end{aligned}
\end{equation}

In order to make use of gradient-based optimization methods, we will need to make further assumptions about the nature of $\ell$ and $f$ that will allow us to compute a gradient. The first assumption that we make ensures that $x(t)$ is continuously differentiable with respect to $\theta$.
\begin{assumption} \label{assumption-regularity-dynamics}
Suppose that the values of $x(t)$ are restricted to a domain, $\mathcal{D}$; that $f$ is continuous in $t$ and uniformly Lipschitz continuous on $\mathcal{D}$; and, finally, that $f$ has a derivative with respect to its second argument, denoted by $f_x$, and $f_x$ is continuous in both of its arguments on $[0,T]\times \mathcal{D}$.
\end{assumption}
Then, by Theorem 7.2 of \cite{coddington1955theory}, $x(t)$ is continuously differentiable with respect to its initial value that takes value in $\mathcal{D}$, for all $t \in [0,T]$. We denote the derivative of $x$ with respect to $\theta$ by $x_\theta$. Moreover, by (7.13) and (7.14) of \cite{coddington1955theory}, the following is valid:
\begin{equation} \label{eqn-forward-derivative-dynamics}
\left\lbrace
\begin{aligned}
\dot x_\theta (t) &= f_x(t, x(t) ) x_\theta (t), ~ t \in [0,T] \\
x_\theta(0) &= I_q,
\end{aligned}
\right.
\end{equation}
where $I_q$ is the $q \times q$ identity matrix.

We will also need to make further assumptions about the nature of $\ell$ in order to compute the gradient of \cref{eqn-optimization-problem}. These assumptions will be more general for the stochastic gradient descent method in comparison to those for the Kalman-based Stochastic Gradient Descent method. We will introduce these assumptions in each of the subsections below, where they are most relevant.

\subsection{Stochastic Gradient Descent for Data Assimilation}

In order to compute the stochastic gradient, we will begin by first computing the gradient. In turn, in order to compute the gradient, we will need a differentiability assumption of $\ell$.
\begin{assumption} \label{assumption-loss-differentiable}
Suppose that $\ell$ is continuously differentiable with respect to its second argument, which we denote by $\ell_x$. 
\end{assumption}

Then the gradient of $G(\theta)$ in \cref{eqn-optimization-problem}, which we denote by $\nabla G$ is supplied by
\begin{equation} \label{eqn-gradient-forward}
\nabla G(\theta) = \sum_{i=1}^N x_\theta(t_i)' \ell_x(y_i, x(t_i)),
\end{equation}
satisfying \cref{eqn-dynamics-ivp,eqn-forward-derivative-dynamics}. While it is possible to simultaneously integrate \cref{eqn-dynamics-ivp,eqn-forward-derivative-dynamics} to compute \cref{eqn-gradient-forward}, an efficient implementation would require integrating $d+1$, $p$-dimensional vectors. And while this may be justifiable in some contexts, an alternative and mathematically equivalent approach is to use the adjoint system \cite{griewank2012invented}, which is given by
\begin{equation} \label{eqn-adjoint-dynamics}
\left\lbrace
\begin{aligned}
\dot \chi (t) &= - f_x(t, x(t) )' \chi(t), ~ t \in [0,T] \\
	 \chi(t_i^-) &= \chi(t_i^+) + \ell_x(y_i, x(t_i)), ~ i=1,\ldots,N \\
	 \chi (T) &= 0,
\end{aligned}
\right.
\end{equation}
where $t_i^-$ is the left-hand limit to $t_i$ and $t_i^+$ is the right-hand limit to $t_i$. Specifically, by first integrating \cref{eqn-dynamics-ivp} and then integrating \cref{eqn-adjoint-dynamics} \textit{backwards} in time, we have that $\nabla G(\theta) = \chi(0)$. Note that when the numerical integral is calculated, there is still a substantial storage burden, which can be partially alleviated by checkpointing and re-integration \cite{griewank2012invented}.

Whether the forward method is utilized (i.e., integrating \cref{eqn-dynamics-ivp,eqn-forward-derivative-dynamics}) or the backward method is utilized (i.e., integrating \cref{eqn-dynamics-ivp,eqn-adjoint-dynamics}) will depend on the problem dimensions, hardware and practical constraints. Regardless of the method, an unbiased stochastic gradient can be calculated in both cases as we now show.

\begin{lemma} \label{lemma-stoc-grad-unbiased}
Suppose \cref{assumption-regularity-dynamics,assumption-loss-differentiable} hold.
Let $\mathcal{S}$ be a randomly generated subset of $\lbrace 1,\ldots,N \rbrace$. For any $s \in \lbrace 1,\ldots, N \rbrace$, let $\pi_s$ denote the probability that $s \in \mathcal{S}$. Then, an unbiased stochastic gradient, denoted $\nabla g_\mathcal{S}(\theta)$, using the forward method is supplied by
\begin{equation}
\nabla g_{\mathcal{S}}(\theta) = \sum_{s \in \mathcal{S}} \frac{1}{\pi_s} x_\theta(t_s)' \ell (y_s, x(t_s) ),
\end{equation}
satisfying \cref{eqn-dynamics-ivp,eqn-forward-derivative-dynamics}. Equivalently, $g_\mathcal{S}(\theta)$ can be computed by a backwards integration of
\begin{equation}
\left\lbrace
\begin{aligned}
\dot \chi (t) &= - f_x(t, x(t) )' \chi(t), ~ t \in [0,T] \\
	 \chi(t_s^-) &= \chi(t_s^+) + \frac{1}{\pi_s} \ell_x(y_s, x(t_s)), ~ s \in \mathcal{S} \\
	 \chi (T) &= 0,
\end{aligned}
\right.
\end{equation}
where first \cref{eqn-dynamics-ivp} is integrated, and setting $g_\mathcal{S}(\theta) = \chi(0)$. 
\end{lemma}
\begin{proof}
It is rather straightforward to show that the backward method produces the same values as the forward method. Therefore, we need only compute the expected value of the forward method to show that the gradient is unbiased, which we can do directly with a simple trick. First, let $\1{\cdot}$ denote the indicator that the specified event occurs (i.e., it is one if the event specified holds, and zero otherwise). Then,
\begin{align}
\E{\nabla g_\mathcal{S}(\theta) } &= \E{ \sum_{s \in \mathcal{S}} \frac{1}{\pi_s} x_\theta(t_s)' \ell_x (y_s, x(t_s) ) } \\
						  &= \E{ \sum_{i=1}^N \frac{1}{\pi_i} x_\theta(t_i)' \ell_x (y_i, x(t_i)) \1{ i \in \mathcal{S} } } \\
						  &= \sum_{i=1}^N \frac{1}{\pi_i} x_\theta(t_i)' \ell_x (y_i, x(t_i)) \Prb{ i \in \mathcal{S} } \\
						  &= \sum_{i=1}^N x_\theta(t_i)' \ell_x (y_i, x(t_i)).
\end{align}
\end{proof}

It is worth remarking here about the inherent differences between computing a stochastic gradient in data assimilation context versus computing a stochastic gradient in the machine learning context, and how these differences require a more thoughtful approach when computing stochastic gradients in the data assimilation context. When comparing these two contexts, the primary difference stems from the cost of evaluating the gradient of the loss function from a given observation. In the machine learning context, the cost of evaluating the gradient of the loss function for a given observation is the sum of the cost of moving the observation from memory (or long-term storage) to cache, and then evaluating the gradient of the loss function. In the data assimilation context, the total cost of evaluating the gradient of a loss function includes the costs in the machine learning context, plus the cost of integrating the state vector, $x$, and its derivative, either $x_\theta$ or $\chi$. In the data assimilation context, this latter cost is the dominant cost owing to the expense of the numerical integration scheme and, perhaps more importantly, the additional storage requirements necessary for calculating the gradients.

For this reason, the sampling process for computing a stochastic gradient in the data assimilation process \textit{cannot} be as arbitrary as it can be for the machine learning context. To understand precisely what we mean, consider the potential costs of computing the data assimilation stochastic gradient using the defacto machine learning standard of a simple random sample of size from a collection of $N$ equally spaced observations over the interval $[0,T]$. Then the time difference between two sequential observations in $\mathcal{S}$, the random sample, can vary from the $\mathcal{O}(T/N)$ to $\mathcal{O}(T)$. In the lower extreme, our numerical integration step size is determined by the separation of observations rather than the stability of the numerical integration scheme, which, in the high-frequency observation regime, is debilitatingly small (i.e., $N$ is very large relative to $T$ and $\mathcal{O}(T/N)$ is much smaller than the preferred integration step size). In the larger extreme, our numerical integration step sizes are (most likely) coinciding with time points where observations have occurred, but these are not being used for the gradient calculation as they do not belong in $\mathcal{S}$; that is, we are paying the price of the numerical integration but ignoring the observations occurring at time points that are not in $\mathcal{S}$. 
Thus, we are either required to make unnecessary numerical integration steps, or we inefficiently omit observations that are can improve the estimation as they do not belong to the sample set, $\mathcal{S}$.

As the above discussion highlights, we want to use sampling schemes that fully exploit the costs of the numerical integration procedure. Therefore, the simplest appropriate sampling scheme is to use systematic random samples (see  Ch. 8 of \cite{cochran2007sampling}). To understand this, consider, again, $N$ evenly spaced observations over the interval $[0,T]$, and suppose that the preferred, application-dependent numerical integration step is given by $\kappa T / N$, 
where $\kappa \in \mathbb{Z}_{>1}$. Then, in the systematic sampling scheme, we randomly choose between the first and $\kappa^\mathrm{th}$ observation, and then choose every $\kappa^{\mathrm{th}}$ observation thereafter. For example, if we were to select the first observation, then the set of observations would have indices in
$\lbrace 1 + z \kappa: z \in \mathbb{Z} \rbrace \cap \lbrace 1,\ldots, N \rbrace.
$
Thus, in the systematic sampling case, we see that each numerical integration step corresponds to an observation time point and the resulting stochastic gradient is unbiased if we use the rescaling in \cref{lemma-stoc-grad-unbiased}.\footnote{There are sampling schemes that are closely related to systematic sampling that can also be used. See \cite{cochran2007sampling}.}

\begin{remark}
In the above example of systematic sampling, we assumed that the integration time steps corresponded to observation time points. In the high-frequency observation regime, this feature would likely hold approximately, and, thus, the reasoning would still apply. 
\end{remark}

Of course, computing the search direction efficiently is only part of the problem. Indeed, computing the step sizes appropriately is also an important consideration. In fact, choosing non-negative step sizes, $\lbrace \eta_k \rbrace$, that satisfy the Robins and Monro conditions \cite{robbins1951stochastic} of
\begin{equation} \label{robbins-monro}
\sum_{k=0}^\infty \eta_k = \infty \quad\text{and}\quad \sum_{k=0}^\infty \eta_k^2 < \infty,
\end{equation}
will lead to convergence under rather general nonconvex structures.\footnote{In our previous work we showed convergence to a stationary point under Bottou-Curtis-Nocedal functions with probability one \cite{patel2020stopping}. If we then take our previous result with that of Theorem 4.9 from \cite{bottou2018optimization}, then one can show that a stationary point is also achieved in $L^2$ by Fatou's lemma.} Unfortunately, these nonconvex structures do \textit{not} sufficiently characterize the data assimilation problem satisfying \cref{assumption-regularity-dynamics}---one can directly see this by trying to find bounds (i.e., (4.8) in \cite{bottou2018optimization} or NM-3 in \cite{patel2020stopping}) on the variance of a stochastic gradient for an arbitrary parameter in the domain. Moreover, the conditions stated in \cref{robbins-monro} must be satisfied carefully, else a poor selection can lead to highly undesirable rates of convergence \cite{nemirovski2009robust} or numerical instability \cite{patel2017impact}. We will not take up the issue of convergence here for the structure imposed by the data assimilation problem, nor the methodology of how to choose the step sizes appropriately, which is a very active area of research. Here, we simply state the procedure in \cref{algorithm-sgd}.

\begin{algorithm}
\caption{SGD for High-frequency Observation Data Assimilation}
\label{algorithm-sgd}
\begin{algorithmic}[1]
\REQUIRE $\theta \in \mathbb{R}^p$
\REQUIRE $\lbrace \eta_k \rbrace$ satisfying \cref{robbins-monro}
\STATE $k \leftarrow 0$
\WHILE{stopping criteria is not satisfied (see \cite{patel2020stopping})}
	\STATE Independently generate $\mathcal{S}$ using a systematic or stratified sampling procedure
	\STATE Compute $\nabla g_\mathcal{S}(\theta)$ using \cref{lemma-stoc-grad-unbiased}
	\STATE $\theta \leftarrow \theta - \eta_k \nabla g_\mathcal{S}(\theta)$
	\STATE $k \leftarrow k + 1$
\ENDWHILE
\RETURN $\theta$
\end{algorithmic}
\end{algorithm}

\subsection{Kalman-based Stochastic Gradient Descent for Data Assimilation}

We now consider second-order stochastic gradient methods, which are known to be highly robust to the choice of step sizes as demonstrated in empirical studies \cite{xu2020second}. One second-order stochastic gradient method, Kalman-based Stochastic Gradient Descent (kSGD), has been shown to hold this property theoretically in a limited context \cite{patel2016kalman}. Owing to this reason and its ability to take advantage of common structures in data assimilation problems (specified below), we discuss how to implement kSGD for data assimilation. 

As a first step, we must guarantee that the second derivative of the dynamics exists even though it will \textit{not} be explicitly calculated for kSGD. Towards this end, the following assumption strengthens \cref{assumption-regularity-dynamics}, and the corresponding differentiability of the solution, $x(t)$, follows from Theorem 7.5 of \cite{coddington1955theory}.

\begin{assumption} \label{assumption-regularity-dynamics-second}
Suppose that the values of $x(t)$ are restricted to a domain, $\mathcal{D}$; that $f$ is continuous in $t$ and uniformly Lipschitz continuous on $\mathcal{D}$; and, finally, that $f$ twice differentiable with respect to its second argument, where the first and second derivatives are denoted by $f_x$ and $f_{xx}$ respectively, and $f_{xx}$ is continuous in both of its arguments on $[0,T]\times \mathcal{D}$.
\end{assumption}

Just as for the first-order stochastic gradient case, we will also need the twice differentiability of the loss function, even though it will \textit{not} be explicitly calculated for kSGD. Moreover, we will restrict ourselves to specifically structured loss functions that are overwhelmingly common in data assimilation. In order to state this special structure, we will explicitly denote the dependence of $x(t)$ and its derivatives on the parameter $\vartheta$ by $x(t;\vartheta)$, $x_\theta(t; \vartheta)$ and $x_{\theta\theta}(t; \vartheta)$. We now state the assumptions that dictate this special structure and follow it with an example that is common in data assimilation. We then discuss how we will leverage this structure.

\begin{assumption} \label{assumption-quasi-likelihood}
Suppose that $\ell$ is twice continuously differentiable with respect to its second argument, for which we denote the first and second derivatives by $\ell_{x}$ and $\ell_{xx}$, respectively. Suppose that, for some $\theta^*$ and for each $i=1,\ldots,N$, $y_i$ has mean $\mu(x(t_i;\theta^*))$, where $\mu$ is a smooth, known function of its argument; and $y_i$ has variance $V( x(t_i;\theta^*) )$, where $V$ is a smooth function of its argument and takes positive definite matrix values. Finally, suppose that
\begin{equation} \label{eqn-quadlikelihood}
\ell_{x}(y_i, x(t_i; \theta^*)) = -\mu_x(x(t_i;\theta^*))'V(x(t_i;\theta^*))^{-1}\left[y_i - \mu(x(t_i;\theta^*)\right],
\end{equation}
for $i = 1,\ldots,N$, where $\mu_x$ is the derivative of $\mu$.
\end{assumption}

A plurality of loss functions used in data assimilation satisfy \cref{assumption-quasi-likelihood}, and fall under the umbrella of Generalized Estimating Equations \cite{hardin2003generalized}. For a simple example of how this assumption is satisfied, we consider the case of the observations being normally distributed about a nonlinear transformation of the state variables.

\begin{example}
Let $\lbrace (t_i, y_i) : i = 1,\ldots,N \rbrace$ be observations such that $y_i$ are normally distributed with mean $\mu( x(t_i; \theta^*) )$ for some known, smooth function $\mu$ and specific parameter $\theta^*$, and variance $V$ which is invertible. If $\ell(y_i,x(t_i;\theta^*))$ is the negative loglikelihood corresponding to the Gaussian distribution, then
\begin{equation}
\ell_x (y_i, x(t_i; \theta^*) ) = -\mu_x ( x(t_i; \theta^*) )'V^{-1}\left[ y_i - \mu( x(t_i; \theta^*) ) \right],
\end{equation}
where $\mu_x$ denotes the derivative of $\mu$ with respect to its argument. 
\end{example}

We now turn to how \cref{assumption-quasi-likelihood} is used. \cref{assumption-quasi-likelihood} is a foundational structural assumption exploited in Gauss-Newton algorithms (e.g., see pages 246 to 247 of \cite{nocedal2006numerical}) owing to the following fact \cite{wedderburn1974quasi}: for any $\ell$ and observations satisfying \cref{assumption-quasi-likelihood} for a given $\theta^*$,
\begin{equation} \label{eqn-quasi-hessian}
\begin{aligned}
\E{ \ell_x(y_i, x(t_i;\theta^*) \ell_x (y_i, x(t_i;\theta^*))' } &= - \E{ \ell_{xx}(y_i, x(t_i;\theta^*) } \\ &= \mu_x( x(t_i;\theta^*) )' V( x(t_i;\theta^*) )^{-1} \mu_x ( x(t_i; \theta^* ) ),
\end{aligned}
\end{equation}
for $i = 1,\ldots,N$. That is, the expected Hessian of the loss at the generating parameter, $\theta^*$, is characterized through the gradient of the loss function. 

To see how this structure is exploited in kSGD, let $\mathcal{S}$ be a random subset of $\lbrace 1,\ldots,N \rbrace$, and let 
\begin{equation}
g_\mathcal{S} (\theta) = \sum_{i \in \mathcal{S}} \frac{1}{\pi_i} l(y_i, x(t_i;\theta) ),
\end{equation}
where the selection of $\mathcal{S}$ is subject to the same considerations as the SGD case; and we recall that $\pi_i = \Prb{ i \in \mathcal{S}}$. Following the strategy outlined in \S 4.1 of \cite{patel2018identification}, a parameter estimate $\theta_{k-1}$ with a corresponding symmetric, positive definite matrix $C_{k-1} \in \mathbb{R}^{q \times q}$ can be updated using the information in $g_\mathcal{S} (\theta)$ by solving
\begin{equation}
\begin{aligned}
\min_{\theta} &~g_\mathcal{S}(\theta) + \frac{1}{2} \norm{ \theta - \theta_{k-1}}_{ C_{k-1}^{-1}}^2,
\end{aligned}
\end{equation}
subject to \cref{eqn-dynamics-ivp}. However, this optimization problem is usually as hard as the original problem, and can be simplified by approximating $g_\mathcal{S}$ with its second-order Taylor expansion and then using \cref{eqn-quasi-hessian} in place of the Hessian of $g_{\mathcal{S}}$. To show this, let $i_1,\ldots,i_{|\mathcal{S}|}$ denote the indices in $\mathcal{S}$ in increasing order, and let
\begin{equation}
r_{\mathcal{S}}(\theta) = \begin{bmatrix}
y_{i_1} - \mu( x(t_{i_1};\theta) ) \\
\vdots  \\
y_{i_{|\mathcal{S}|}} - \mu ( x(t_{i_{|\mathcal{S}|}}; \theta) ) 
\end{bmatrix}, \quad
D_{\mathcal{S}}(\theta) = \begin{bmatrix}
\mu_x( x(t_{i_1};\theta) ) x_\theta( t_{i_1}; \theta) \\
\vdots \\
\mu_x( x(t_{i_{|\mathcal{S}|}};\theta) ) x_\theta( t_{i_{|\mathcal{S}|}}; \theta)
\end{bmatrix},
\end{equation}
and
\begin{equation}
W_{\mathcal{S}}^{-1}(\theta) = \begin{bmatrix}
\frac{1}{\pi_{i_1}} V( x(t_{i_1};\theta ) )^{-1} & & 0 \\
&  \ddots & \\
0 & &  \frac{1}{\pi_{i_{|\mathcal{S}|}}} V ( x(t_{i_{|\mathcal{S}|}} ;\theta) )^{-1}
\end{bmatrix}.
\end{equation}

Then, when we replace $g_{\mathcal{S}}(\theta)$ with the stated approximation, the minimization problem reduces to
\begin{equation}
\begin{aligned}
\min_{\theta} ~ &- r_{\mathcal{S}}(\theta_{k-1})' W_{\mathcal{S}}^{-1}(\theta_{k-1}) D_{\mathcal{S}}(\theta_{k-1})(\theta - \theta_{k-1})  \\
			    & + \frac{1}{2}( \theta - \theta_{k-1} )' \left[ D_{\mathcal{S}}(\theta_{k-1})' W_{\mathcal{S}}^{-1}(\theta_{k-1}) D_{\mathcal{S}}(\theta_{k-1}) + C_{k-1}^{-1} \right] (\theta - \theta_{k-1}),
\end{aligned}
\end{equation}
which can be explicitly solved to give the update
\begin{equation} \label{eqn-ksgd-update-1}
\theta_{k} = \theta_{k-1} + \left[ D_{\mathcal{S}}(\theta_{k-1})' W_{\mathcal{S}}^{-1}(\theta_{k-1}) D_{\mathcal{S}}(\theta_{k-1}) + C_{k-1}^{-1} \right]^{-1} D_{\mathcal{S}}(\theta_{k-1})'W_{\mathcal{S}}^{-1}(\theta_{k-1}) r_{\mathcal{S}}(\theta_{k-1}).
\end{equation}
With an application of the Woodbury matrix identity, we also recover the equivalent form,
\begin{equation} \label{eqn-ksgd-update-2}
\theta_{k} = \theta_{k-1} + C_{k-1} D_{\mathcal{S}}(\theta_{k-1})'\left[ W_{\mathcal{S}}(\theta_{k-1}) + D_{\mathcal{S}}(\theta_{k-1}) C_{k-1} D_{\mathcal{S}}(\theta_{k-1})' \right]^{-1} r_{\mathcal{S}}(\theta_{k-1}).
\end{equation}
Moreover, we will then compute $C_{k}$ from $C_{k-1}$ by
\begin{equation} \label{eqn-ksgd-precision-update}
C_{k}^{-1} = D_{\mathcal{S}}(\theta_{k-1})' W_{\mathcal{S}}^{-1}(\theta_{k-1}) D_{\mathcal{S}}(\theta_{k-1}) + C_{k-1}^{-1}.
\end{equation}
When we choose $\theta_0$ and $C_0 = I$, we have a complete procedure specified by either \cref{eqn-ksgd-update-1,eqn-ksgd-precision-update} or \cref{eqn-ksgd-update-2,eqn-ksgd-precision-update}. We summarize the procedure in \cref{algorithm-ksgd}.

\begin{algorithm}
\caption{kSGD for High-frequency Observation Data Assimilation}
\label{algorithm-ksgd}
\begin{algorithmic}[1]
\REQUIRE $\theta \in \mathbb{R}^p$
\STATE $k \leftarrow 0$
\STATE $C \leftarrow I$
\WHILE{stopping criteria is not satisfied (see \cite{patel2020stopping})}
	\STATE Independently generate $\mathcal{S}$ using a systematic or stratified sampling procedure
	\STATE Compute $r_\mathcal{S}$, $W_{\mathcal{S}}$ or $D_{\mathcal{S}}$ using \cref{eqn-dynamics-ivp,eqn-forward-derivative-dynamics}
	\STATE Update $\theta$ by \cref{eqn-ksgd-update-1} or \cref{eqn-ksgd-update-2}
	\STATE Update $C^{-1}$ by \cref{eqn-ksgd-precision-update}
	\STATE $k \leftarrow k + 1$
\ENDWHILE
\RETURN $\theta$
\end{algorithmic}
\end{algorithm}

The two procedures listed above are mathematically equivalent, but the choice between them is often dependent on computational considerations. For example, in the high-frequency observation domain, we might consider $|\mathcal{S}|$ to be reasonably large, and, consequently, this would result in the dimension of $r_{\mathcal{S}}(\theta_{k-1})$ to be large. In this context, the implicit linear system in \cref{eqn-ksgd-update-1} is most likely substantially smaller than that in \cref{eqn-ksgd-update-2}. When the dimension of the parameter is also large, we will need to develop techniques that can approximate the large matrices that will result and compute the updates (e.g., \cite{petra2019structured}).

As it may now be apparent, kSGD is rather closely related to adaptations of the Kalman Filter for parameter estimation, namely the augmented Kalman Filter (e.g., see \cite{ljung1977analysis}) and the incremental optimizer studied in \cite{bertsekas1996incremental}. Therefore, we conclude this section by briefly contrasting these to kSGD. The augmented Kalman filter augments the state vector with the parameter vector, and filters over both of these quantities. However, this would only be an approximate solution to the variational assimilation problem, whereas kSGD attempts to solve the variational problem.\footnote{We do will show this numerically in this work rather than theoretically, which we leave to future efforts.} The incremental optimizer in \cite{bertsekas1996incremental} is a deterministic special case of kSGD that includes a forgetting parameter which can also be incorporated here. However, the analysis in \cite{bertsekas1996incremental} relies on a notion of persistent excitation, which is often inappropriate in data problems and would motivate an alternative analysis (e.g., \cite{patel2016kalman}) that we leave to future work.

%% file: sections/Experiments.tex
Recall, high-frequency observations---that is, observations that occur at regular intervals smaller than a desired numerical integration step---result in substantial computational requirements when used in data assimilation tasks. While the computational burden induced by high-frequency observations can be addressed using data modification strategies (\cref{subsection-data-modification}), these data-modification strategies will often result in a non-trivial reduction in the quality of the estimates (\cref{subsection-data-modification-consequences}). Motivated by this apparent computational-statistical trade-off, we adapted stochastic approximation methods---a first-order method (Stochastic Gradient Descent) and a second-order method method (Kalman-based Stochastic Gradient Descent)---to the data assimilation problem with high-frequency observations that we claimed would circumvent this trade-off: that is, these methods would respect computational budgets while also achieving high statistical accuracy. 

In this section, we provide numerical evidence to show that, indeed, our adaptations of these two stochastic approximation methods to the data assimilation problem will achieve performance comparable to, or better than, a deterministic method on the unmodified optimization problem while also respecting computational budgets. To this end, we describe our general experimental setup and our test dynamical systems in \cref{subsec-experimental-setup}. Then, in \cref{subsec-comparison}, we compare first-order methods and we compare second-order methods.

\subsection{Experimental Setup}
\label{subsec-experimental-setup}

We consider three dynamical systems models: the FitzHugh-Nagumo relaxation oscillator \cite{rocsoreanu2012fitzhugh}, Lotka-Volterra predator-prey model \cite{lotka1920analytical}, and Van der Pol relaxation oscillator \cite{van1960theory}. Specifically, we investigate FitzHugh-Nagumo model on time interval [0,50] with integration step size of one time unit and observations occurring every 0.01 time units; we implement Lotka-Volterra model on time interval [0,10] with integration step size 0.5 time units and observations occuring every 0.005 time units; we implement Van der Pol model on time interval [0,10] with integration step size 0.1 time units and observations occurring every 0.001 time units. For each model, we choose the ``true'' parameter, the initial state, and integration interval based on the literature; we choose a preferred, large, stable numerical integration step size for Ralston's fourth order explicit Runge-Kutta method by numerical comparisons with smaller step sizes; and we simulate observations from a linear transformation of the state with additive noise at time points that are one hundredth of the preferred step size. 

We then try to identify the parameter of the model by a collection of first-order methods and second-order methods. The first order methods include our adaptation of Stochastic Gradient Descent (SGD), Gradient Descent (GD) applied to the problem with the aforementioned six data modification schemes (effectively looking only at one percent of the data), and GD applied to the original problem. The second order methods include our adaption of Kalman-based Stochastic Gradient Descent (kSGD), Gauss-Newton (GN) applied to the original problem with the aforementioned six data modification scheme, and GN applied to the original problem. 

To point out some specific salient features of the experiment, we use the forward method to compute the derivatives as using the backward method would require a greater memory burden owing to the large number of observations (even with checkpointing) because of the size of the problem. Additionally, the step-sizes for each of the above methods were hand-tuned and non-adaptive: this was to ensure competitiveness across these methods as adaptively selecting such parameters incurs a substantial computational burden for deterministic methods, and such methods are still being actively developed for stochastic methods.

As described previously, we now consider the second route of experimentation as discussed in \cref{subsection-data-modification-consequences}. That is, we allow the solvers to run for a fixed time budget (about one second based on the size of the problems) and compare the progress of the methods on this fixed computational budget. For reporting and comparison purposes, we plot the relative error defined in \cref{eqn-relative-accuracy} in comparison to elapsed compute time. 

\subsection{Comparison of First-order and Second-order Methods}
\label{subsec-comparison}

We will begin by first comparing within first-order methods, followed by a comparison within second-order methods, and then brief comments comparing across the different orders.

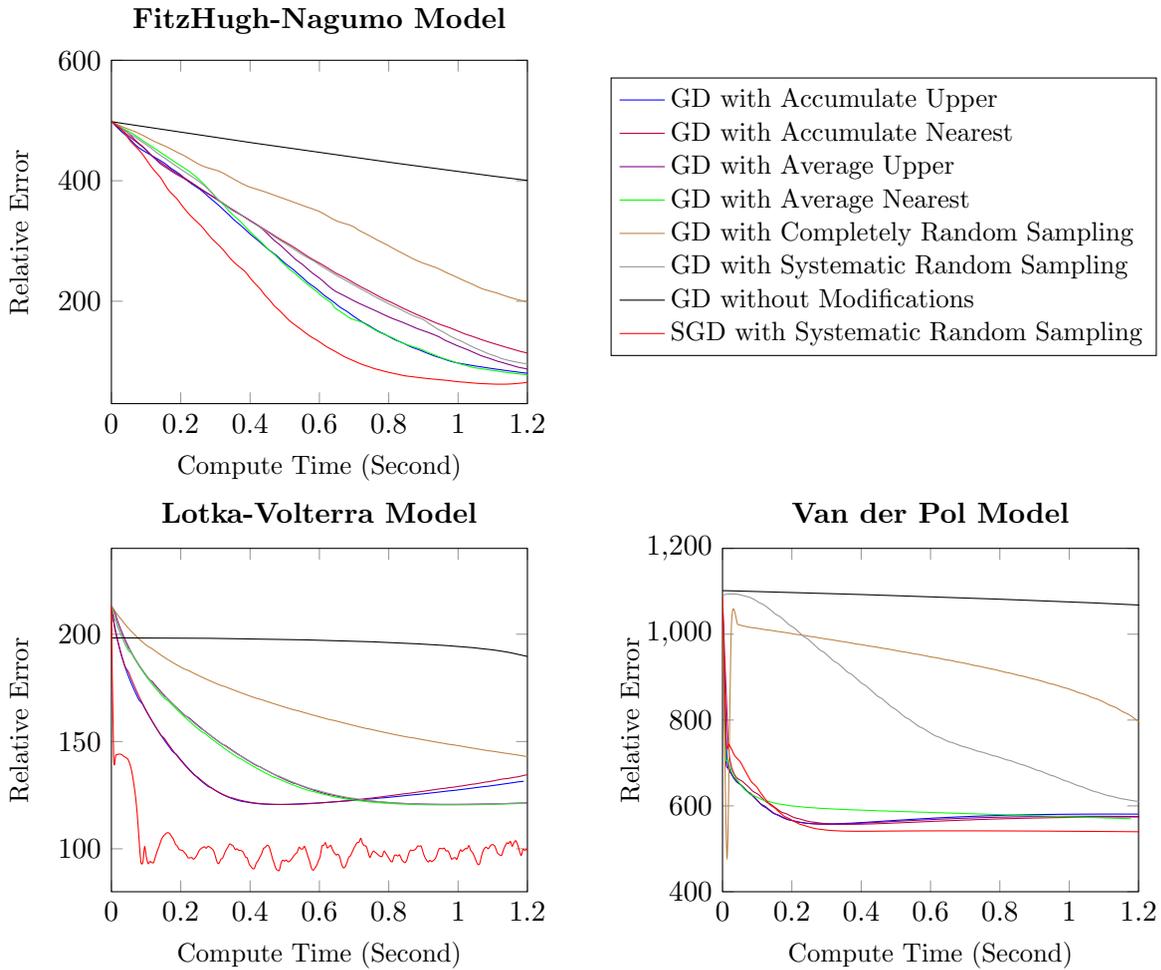
\begin{figure}[h]
	\centering
	\begin{subfigure}[b]{\linewidth}
		\input{tikz/FN_SGD}
		\label{figure:FN_SGD}
	\end{subfigure}
	
	\begin{subfigure}{0.48\linewidth}
		\input{tikz/LV_SGD}
		\label{figure:LV_SGD}
	\end{subfigure}
	\hfill
	\begin{subfigure}{0.48\linewidth}
		\input{tikz/vanderpol_SGD}
		\label{figure:vanderpol_SGD}
	\end{subfigure}
	
	\caption{A comparison of first-order optimizers and optimization problems on the three test models.}
	\label{figure:first-order-methods}
\end{figure} 

As shown in \cref{figure:first-order-methods}, SGD with systematic random sampling attains the smallest relative error at almost every time point to various degrees for the three dynamics, while the least progress is made by GD on the unmodified problem. Interestingly, GD on many of the data modified problems tends to perform rather comparably to SGD, but, as we saw in \cref{subsection-data-modification-consequences}, may have a higher limit to their lowest achievable error if a longer budget was allowed. It is also worth pointing out that these estimation problems are nonconvex, and, in the case of SGD, we observed that the iterates are converging to local minima. In practice, this nonconvexity is dealt with in two ways. One way is by choosing an initialization that is based on previous calculations or the physics of the system---here, we chose an uninformed initialization by not taking into account any prior information. The other way of handling the nonconvexity is to use adaptive step sizes that ensure a certain amount of reduction in the objective function is achieved. Of course, adaptive step sizes for stochastic methods are still an active area of research, and will not be something that we address here. 

\begin{figure}[ht]
	\centering
	\begin{subfigure}[b]{\linewidth}
		\input{tikz/FN_kSGD}
		\label{figure:FN_kSGD}
	\end{subfigure}
	\begin{subfigure}{0.48\linewidth}
		\input{tikz/LV_kSGD}
		\label{figure:LV_kSGD}
	\end{subfigure}
	\hfill
	\begin{subfigure}{0.48\linewidth}
		\input{tikz/vanderpol_kSGD}
		\label{figure:vanderpol_kSGD}
	\end{subfigure}
	
	\caption{A comparison of second-order optimizers and optimization problems on the three test models.}
	\label{figure:second-order-methods}
\end{figure}
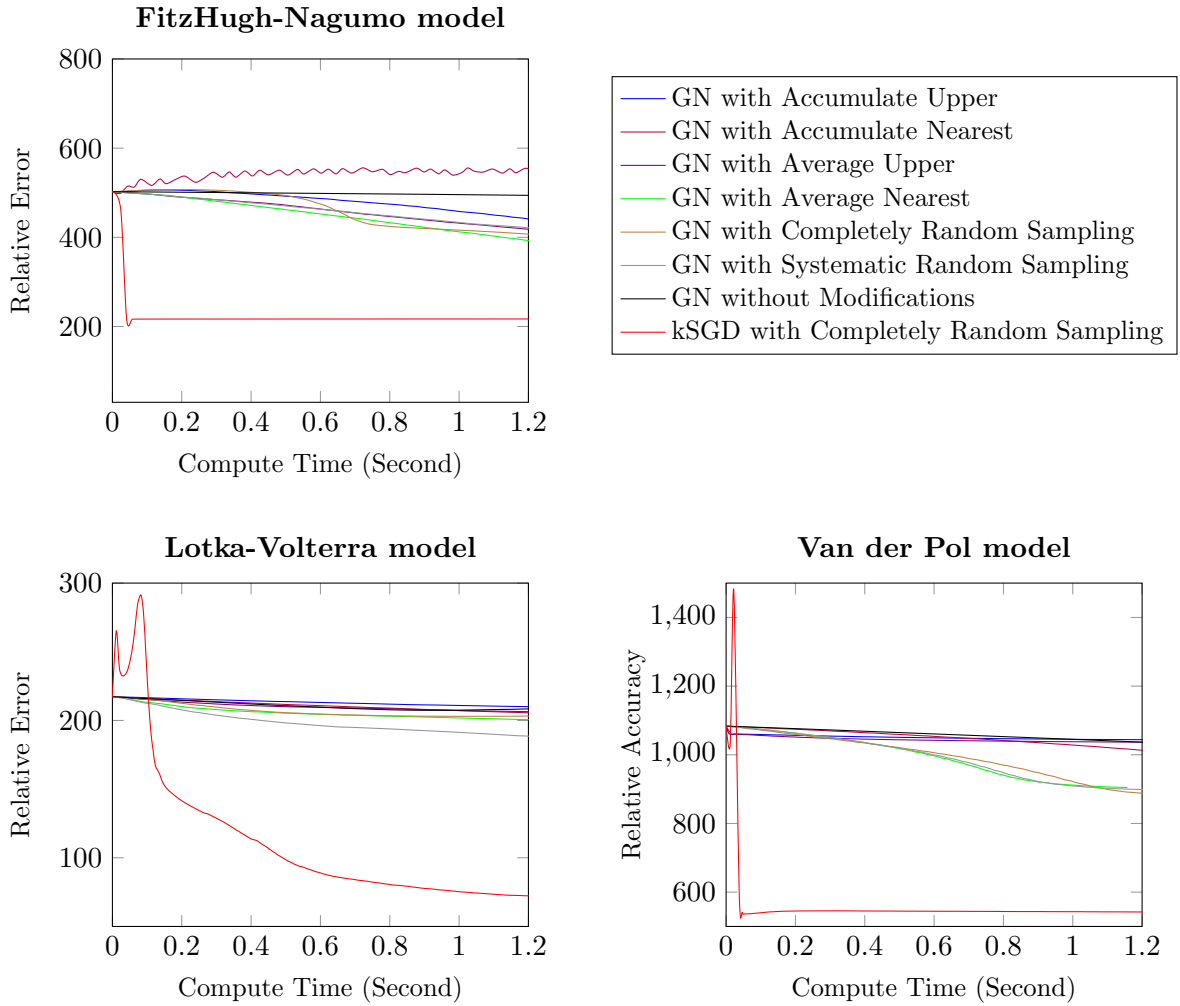 

A similar story for second-order methods is shown in \cref{figure:second-order-methods}. Specifically, kSGD with systematic random sampling attains the best statistical accuracy among all the second order methods for all three dynamics within the limited computational budgets. In contrast to the first-order results, GN on the data modified problems do not perform comparably to kSGD. We can attribute this underpeformance of GN on the modified data problems to the increased computational costs associated with second-order methods---which are usually only justified as the iterates approach the optimal point---, and to potentially over-optimizing to the modified data objective function as we saw in \cref{subsection-data-modification-consequences}. Again, the considerations around nonconvexity also play a role for these second-order methods. If we use the first route mentioned in \cref{subsection-data-modification-consequences}, we will be able to show that SGD and kSGD circumvent the statistical-computational trade-off, which will be left for review. 

When we compare across first-order and second-order methods, we notice that kSGD seems to find minima faster than SGD (with the exception of the Lotka-Volterra problem). This seems to suggest that if kSGD can be feasibly implemented then it seems to outperform its first-order counterpart on these problems. The downside, however, is that there are more parameters to tune of kSGD, which reinforces the need for adaptive techniques for determining these parameters. Again, we will leave such efforts to future work.

To summarize, we indeed confirm that stochastic approximation methods appear to offer a way to make substantial progress on data assimilation problems with high-frequency observations while respecting computational budgets in comparison to solving data modified variants of the data assimilation problem.

%% file: tikz/FN_SGD.tex
\begin{tikzpicture}
			\begin{axis}[
				title={\bf FitzHugh-Nagumo Model},
				xlabel={\small Compute Time (Second)},
				ylabel={\small Relative Error},
				table/col sep=comma,
				xmin=0, xmax=1.2,
				ymin=30, ymax=600,
				domain=0:3,
				restrict y to domain=0:600,
				width = 2.8in,
				legend style={at={(1.2,0.95)},anchor=north west,font=\small},
				legend cell align={left}
				]
				
				\addplot[smooth,color=blue] table [x=time,y=error]{"data/FitzHugh_Nagumo_First_Order/GD_allToUpper.csv"};
				\addlegendentry {GD with Accumulate Upper};
				\addplot[smooth,color=purple] table [x=time,y=error]{"data/FitzHugh_Nagumo_First_Order/GD_allToNearest.csv"};
				\addlegendentry {GD with Accumulate Nearest};
				\addplot[smooth,color=violet] table [x=time,y=error]{"data/FitzHugh_Nagumo_First_Order/GD_averageToUpper.csv"};
				\addlegendentry {GD with Average Upper};
				\addplot[smooth,color=green] table [x=time,y=error]{"data/FitzHugh_Nagumo_First_Order/GD_averageToNearest.csv"};
				\addlegendentry {GD with Average Nearest};
				\addplot[smooth,color=brown] table [x=time,y=error]{"data/FitzHugh_Nagumo_First_Order/GD_fullRandomSample.csv"};
				\addlegendentry {GD with Completely Random Sampling};
				\addplot[smooth,color=gray] table [x=time,y=error]{"data/FitzHugh_Nagumo_First_Order/GD_systematicRandomSample.csv"};
				\addlegendentry {GD with Systematic Random Sampling};
				\addplot[smooth,color=black] table [x=time,y=error]{"data/FitzHugh_Nagumo_First_Order/GD_nomod.csv"};
				\addlegendentry {GD without Modifications};
				\addplot[smooth,color=red] table [x=time,y=error]{"data/FitzHugh_Nagumo_First_Order/SGD.csv"};
				\addlegendentry {SGD with Systematic Random Sampling};
				
			\end{axis}
\end{tikzpicture}

%% file: tikz/LV_SGD.tex
\begin{tikzpicture}
\begin{axis}[
title={\bf Lotka-Volterra Model},
xlabel={\small Compute Time (Second)},
ylabel={\small Relative Error},
table/col sep=comma,
xmin=0, xmax=1.2,
ymin=80, ymax=240,
domain=0:3,
restrict y to domain=0:600,
width = 2.8in
]

\addplot[smooth,color=blue] table [x=time,y=error]{"data/lotka_volterra_First_Order/GD_allToUpper.csv"};
\addlegendentry {GD with Accumulate Upper};
\addplot[smooth,color=purple] table [x=time,y=error]{"data/lotka_volterra_First_Order/GD_allToNearest.csv"};
\addlegendentry {GD with Accumulate Nearest};
\addplot[smooth,color=violet] table [x=time,y=error]{"data/lotka_volterra_First_Order/GD_averageToUpper.csv"};
\addlegendentry {GD with Average Upper};
\addplot[smooth,color=green] table [x=time,y=error]{"data/lotka_volterra_First_Order/GD_averageToNearest.csv"};
\addlegendentry {GD with Average Nearest};
\addplot[smooth,color=brown] table [x=time,y=error]{"data/lotka_volterra_First_Order/GD_fullRandomSample.csv"};
\addlegendentry {GD with Completely Random Sampling};
\addplot[smooth,color=gray] table [x=time,y=error]{"data/lotka_volterra_First_Order/GD_systematicRandomSample.csv"};
\addlegendentry {GD with Systematic Random Sampling};
\addplot[smooth,color=black] table [x=time,y=error]{"data/lotka_volterra_First_Order/GD_nomod.csv"};
\addlegendentry {GD without Modifications};
\addplot[smooth,color=red] table [x=time,y=error]{"data/lotka_volterra_First_Order/SGD.csv"};
\addlegendentry {SGD with Systematic Random Sampling};

\legend{}; 
\end{axis}
\end{tikzpicture}

%% file: tikz/vanderpol_SGD.tex
\begin{tikzpicture}
\begin{axis}[
title={\bf Van der Pol Model},
xlabel={\small Compute Time (Second)},
ylabel={\small Relative Error},
table/col sep=comma,
xmin=0, xmax=1.2,
ymin=400, ymax=1200,
domain=0:2,
restrict y to domain=0:1200,
width = 2.8in
]

\addplot[smooth,color=blue] table [x=time,y=error]{"data/vanderpol_First_Order/GD_allToUpper.csv"};
\addlegendentry {GD with Accumulate Upper};
\addplot[smooth,color=purple] table [x=time,y=error]{"data/vanderpol_First_Order/GD_allToNearest.csv"};
\addlegendentry {GD with Accumulate Nearest};
\addplot[smooth,color=violet] table [x=time,y=error]{"data/vanderpol_First_Order/GD_averageToUpper.csv"};
\addlegendentry {GD with Average Upper};
\addplot[smooth,color=green] table [x=time,y=error]{"data/vanderpol_First_Order/GD_averageToNearest.csv"};
\addlegendentry {GD with Average Nearest};
\addplot[smooth,color=brown] table [x=time,y=error]{"data/vanderpol_First_Order/GD_fullRandomSample.csv"};
				\addlegendentry {GD with Completely Random Sampling};
				\addplot[smooth,color=gray] table [x=time,y=error]{"data/vanderpol_First_Order/GD_systematicRandomSample.csv"};
				\addlegendentry {GD with Systematic Random Sampling};
				\addplot[smooth,color=black] table [x=time,y=error]{"data/vanderpol_First_Order/GD_nomod.csv"};
				\addlegendentry {GD without Modifications};
				\addplot[smooth,color=red] table [x=time,y=error]{"data/vanderpol_First_Order/SGD.csv"};
				\addlegendentry {SGD with Completely Random Sampling};
				
				\legend{}; 
				
			\end{axis}
\end{tikzpicture}

%% file: tikz/FN_kSGD.tex
\begin{tikzpicture}
\begin{axis}[
title={\bf FitzHugh-Nagumo model},
xlabel={\small Compute Time (Second)},
ylabel={\small Relative Error},
table/col sep=comma,
xmin=0, xmax=1.2,
ymin=30, ymax=800,
domain=0:3,
restrict y to domain=0:800,
width=2.8in,
legend style={at={(1.2,0.95)},anchor=north west,font=\small},
legend cell align={left}
]

\addplot[smooth,color=blue] table [x=time,y=error]{"data/FitzHugh_Nagumo_Second_Order/GN_allToUpper.csv"};
\addlegendentry {GN with Accumulate Upper};
\addplot[smooth,color=purple] table [x=time,y=error]{"data/FitzHugh_Nagumo_Second_Order/GN_allToNearest.csv"};
\addlegendentry {GN with Accumulate Nearest};
\addplot[smooth,color=violet] table [x=time,y=error]{"data/FitzHugh_Nagumo_Second_Order/GN_averageToUpper.csv"};
\addlegendentry {GN with Average Upper};
\addplot[smooth,color=green] table [x=time,y=error]{"data/FitzHugh_Nagumo_Second_Order/GN_averageToNearest.csv"};
\addlegendentry {GN with Average Nearest};
\addplot[smooth,color=brown] table [x=time,y=error]{"data/FitzHugh_Nagumo_Second_Order/GN_fullRandomSample.csv"};
\addlegendentry {GN with Completely Random Sampling};
\addplot[smooth,color=gray] table [x=time,y=error]{"data/FitzHugh_Nagumo_Second_Order/GN_systematicRandomSample.csv"};
\addlegendentry {GN with Systematic Random Sampling};
\addplot[smooth,color=black] table [x=time,y=error]{"data/FitzHugh_Nagumo_Second_Order/GN_nomod.csv"};
\addlegendentry {GN without Modifications};
\addplot[smooth,color=red] table [x=time,y=error]{"data/FitzHugh_Nagumo_Second_Order/kSGD.csv"};
\addlegendentry {kSGD with Completely Random Sampling};

\end{axis}
\end{tikzpicture}

%% file: tikz/LV_kSGD.tex
\begin{tikzpicture}
\begin{axis}[
title={\bf Lotka-Volterra model},
xlabel={\small Compute Time (Second)},
ylabel={\small Relative Error},
table/col sep=comma,
xmin=0, xmax=1.2,
ymin=50, ymax=300,
domain=0:3,
restrict y to domain=0:800,
width=2.8in
]

\addplot[smooth,color=blue] table [x=time,y=error]{"data/lotka_volterra_Second_Order/GN_allToUpper.csv"};
\addlegendentry {GN with Accumulate Upper};
\addplot[smooth,color=purple] table [x=time,y=error]{"data/lotka_volterra_Second_Order/GN_allToNearest.csv"};
\addlegendentry {GN with Accumulate Nearest};
\addplot[smooth,color=violet] table [x=time,y=error]{"data/lotka_volterra_Second_Order/GN_averageToUpper.csv"};
\addlegendentry {GN with Average Upper};
\addplot[smooth,color=green] table [x=time,y=error]{"data/lotka_volterra_Second_Order/GN_averageToNearest.csv"};
\addlegendentry {GN with Average Nearest};
\addplot[smooth,color=brown] table [x=time,y=error]{"data/lotka_volterra_Second_Order/GN_fullRandomSample.csv"};
\addlegendentry {GN with Completely Random Sampling};
\addplot[smooth,color=gray] table [x=time,y=error]{"data/lotka_volterra_Second_Order/GN_systematicRandomSample.csv"};
\addlegendentry {GN with Systematic Random Sampling};
\addplot[smooth,color=black] table [x=time,y=error]{"data/lotka_volterra_Second_Order/GN_nomod.csv"};
\addlegendentry {GN without Modifications};
\addplot[smooth,color=red] table [x=time,y=error]{"data/lotka_volterra_Second_Order/kSGD.csv"};
\addlegendentry {kSGD with Completely Random Sampling};

\legend{};

\end{axis}
\end{tikzpicture}

%% file: tikz/vanderpol_kSGD.tex
\begin{tikzpicture}
\begin{axis}[
title={\bf Van der Pol model},
xlabel={\small Compute Time (Second)},
ylabel={\small Relative Accuracy},
table/col sep=comma,
xmin=0, xmax=1.2,
ymin=500, ymax=1500,
domain=0:3,
restrict y to domain=500:1500,
width=2.8in,
]

\addplot[smooth,color=blue] table [x=time,y=error]{"data/vanderpol_Second_Order/GN_allToUpper.csv"};
\addlegendentry {GN with Accumulate Upper};
\addplot[smooth,color=purple] table [x=time,y=error]{"data/vanderpol_Second_Order/GN_allToNearest.csv"};
\addlegendentry {GN with Accumulate Nearest};
\addplot[smooth,color=violet] table [x=time,y=error]{"data/vanderpol_Second_Order/GN_averageToUpper.csv"};
\addlegendentry {GN with Average Upper};
\addplot[smooth,color=green] table [x=time,y=error]{"data/vanderpol_Second_Order/GN_averageToNearest.csv"};
\addlegendentry {GN with Average Nearest};
\addplot[smooth,color=brown] table [x=time,y=error]{"data/vanderpol_Second_Order/GN_fullRandomSample.csv"};
\addlegendentry {GN with Completely Random Sampling};
\addplot[smooth,color=gray] table [x=time,y=error]{"data/vanderpol_Second_Order/GN_systematicRandomSample.csv"};
\addlegendentry {GN with Systematic Random Sampling};
\addplot[smooth,color=black] table [x=time,y=error]{"data/vanderpol_Second_Order/GN_nomod.csv"};
\addlegendentry {GN without Modifications};
\addplot[smooth,color=red] table [x=time,y=error]{"data/vanderpol_Second_Order/kSGD.csv"};
\addlegendentry {kSGD with Completely Random Sampling};

\legend{};

\end{axis}
\end{tikzpicture}

%% file: sections/Conclusion.tex
When subject to high-frequency observations, solving data assimilation problems is often computationally expensive. Thus, data modification schemes, such as accumulating, averaging and sampling, have often been used to reduce these high computational costs while still supplying some form of parameter and state estimates. Unfortunately, the performance of these traditional modification schemes often sacrifice the quality of the estimates for reductions in computational costs (see \cref{sec:Characterization}).

In an attempt to preserve the statistical quality of the solution while respecting computational budgets, we adapted stochastic gradient descent (SGD) and Kalman-based Stochastic Gradient Descent (kSGD) to the data assimilation problem using systematic random sampling (see \cref{sec:Stochastic}). And, we demonstrated experimentally that our adaptations of SGD and kSGD do indeed provide superior estimates over a fixed computational budget in comparison to solving the data assimilation problem with data modifications (see \cref{sec:Experiments}). Thus, our adaptations of these stochastic approximation methods are able to achieve higher statistical accuracy while respecting computational constraints.

Moving forward, we intend to address several questions that precipitated from this work. First, we will explore how we can efficiently implement adaptive step sizes in the numerical integration step to minimize integration costs and make the most efficient use of the observations. Second, owing to the frequent nonconvexity that we observe for data assimilation problems, we will explore how adaptive step sizes in the outer, numerical optimization subroutines will impact the efficiency of our stochastic approximation methods. Third, we will need to demonstrate the convergence of these stochastic methods to the structures that arise from dynamical systems that are not covered under the current convergence theory of stochastic approximation methods. Finally, we plan to extend our code base to handle larger problems and work with more realistic models and data.